\documentclass{article}
 \usepackage{amsmath}

\newtheorem{thm}{Theorem}[section]
\newtheorem{cor}[thm]{Corollary}
\newtheorem{prop}[thm]{Proposition}

\newtheorem{lem}[thm]{Lemma}

\newtheorem{rem}[thm]{Remark}

\newtheorem{ex}[thm]{Example}

\newcommand{\be}{\begin{equation}}
\newcommand{\ee}{\end{equation}}
\newcommand{\ben}{\begin{enumerate}}
\newcommand{\een}{\end{enumerate}}
\newcommand{\beq}{\begin{eqnarray}}
\newcommand{\eeq}{\end{eqnarray}}
\newcommand{\beqn}{\begin{eqnarray*}}
\newcommand{\eeqn}{\end{eqnarray*}}

\newcommand{\pa}{\partial}

\newcommand{\qed}{\hspace*{\fill}Q.E.D.}  %Use at end of proof

\begin{document}
\title{On a Class of Two-Dimensional Singular Douglas and Projectively flat Finsler
Metrics}
\author{Guojun Yang  }
\date{}
\maketitle
\begin{abstract}
Singular Finsler metrics, such as Kropina metrics and $m$-Kropina
metrics, have a lot of applications in the real world. In this
paper, we study a class of two-dimensional singular Finsler
metrics defined by a Riemann metric $\alpha$ and 1-form $\beta$,
and we characterize those which are Douglasian or locally
projectively flat by some equations. It shows that the main class
induced is an $m$-Kropina metric plus a linear part on $\beta$.
For this class, the local structure of Douglasian or (in part)
projectively flat case is determined, and in particular we show
that a Kropina metric is always Douglasian and a  Douglas
$m$-Kropina metric with $m\ne -1$ is locally Minkowskian. It
indicates that the two-dimensional case is quite different from
the higher dimensional ones.

\


 {\bf MR(2000) subject classification: }
53A20, 53B40
\end{abstract}

\section{Introduction}
       There are
   two important projective invariants in projective Finsler geometry: the Douglas
  curvature ({\bf D}) and the Weyl curvature (${\bf W}^{o}$ in dimension two and
  {\bf W} in higher dimensions) (\cite{Dou}). A Finsler metric is called
  {\it Douglasian} if ${\bf D}=0$. Roughly speaking, a Douglas metric
 is a Finsler metric having the same geodesics as a Riemannian metric.
  A Finsler metric is said to be {\it locally projectively flat} if at every point, there are local
coordinate systems in which geodesics are straight.
  As we know, the locally projectively flat class of Riemannian metrics
  is very limited, nothing but the class
  of constant sectional
  curvature (Beltrami Theorem). However, the class of locally projectively flat Finsler metrics is very
  rich. Douglas metrics form a rich class of Finsler metrics including
locally projectively flat Finsler metrics, and meanwhile there are
many Douglas metrics which are not locally projectively flat.

In this paper, we will concentrate on a special class of
two-dimensional Finsler metrics:  $(\alpha,\beta)$-metrics, and
characterize those which are Douglasian and locally projectively
flat under the condition (\ref{j2}) below. An {\it
$(\alpha,\beta)$-metric} is defined by a Riemannian metric
 $\alpha=\sqrt{a_{ij}(x)y^iy^j}$ and a $1$-form $\beta=b_i(x)y^i$ on a manifold
 $M$, which can be expressed in the following form:
 $$F=\alpha \phi(s),\ \ s=\beta/\alpha,$$
where $\phi(s)$ is a function satisfying certain conditions. It is
known that $F$ is a  regular Finsler metric if $\beta$ satisfies
 $ \|\beta\|_{\alpha} < b_o$ and $\phi(s)$ is  $C^{\infty}$  on $(-b_o,b_o)$ satisfying
 \be\label{j1}
 \phi(s)>0,\quad \phi(s)-s\phi'(s)+(\rho^2-s^2)\phi''(s)>0,
   \quad (|s| \leq \rho <b_o),
 \ee
 where $b_o$ is a positive
constant (\cite{Shen2}). If $\phi(0)$ is not defined or $\phi$
does not satisfy (\ref{j1}),
 then the $(\alpha,\beta)$-metric $F= \alpha \phi
(\beta/\alpha)$ is singular. Singular Finsler metrics have a lot
of applications in the real world (\cite{AIM} \cite{AHM} ). Z.
Shen also introduces singular Finsler metrics in \cite{Shen3}.

Assume that
  $\phi(s)$ is in the following form
   \be\label{j2}
 \phi(s):=cs+s^m\varphi(s),
   \ee
where $c,m$ are constant with $m\ne 0,1$ and $\varphi(s)$ is a
$C^{\infty}$ function on a neighborhood of $s=0$ with
$\varphi(0)=1$, and further for convenience we put $c=0$ if $m$ is
a negative integer. If $m=0$, we have $\phi(0)=1$ and this case
appears in a lot of literatures. When $m\ge 2$ is an integer,
(\ref{j2})
  is equivalent to the following condition
  $$
  \phi(0)=0,\quad
  \phi^{(k)}(0)=0 \ \ (2\le k\le m-1), \quad \phi^{(m)}(0)=m!.
  $$
  Another interesting case is $c=0$ and $\varphi(s)\equiv 1$ in
  (\ref{j2}), and in this case, $F=\alpha\phi(s)$ is called an
  $m$-Kropina metric, and in particular a Kropina metric when $m=-1$.

The case $\phi(0)=1$ has been studied in  a lot of interesting
research papers  (\cite{LSS}--\cite{LS2}, \cite{Ma} \cite{Shen1}
\cite{Shen2} \cite{Y1}--\cite{Y2}). In \cite{LSS} \cite{Shen1},
the authors respectively study and characterize Douglas
$(\alpha,\beta)$-metrics and locally projectively flat
$(\alpha,\beta)$-metrics in dimensions $n\geq 3$ and $\phi(0)=1$,
and further,  the present author solves the case $n=2$ and shows
that the two-dimensional case is quite different from the higher
dimensional ones (\cite{Y1}). In singular case, there are some
papers on the studies of $m$-Kropina metrics and Kropina metrics
(\cite{RR} \cite{ShY} \cite{YN} \cite{YO}). Further, in \cite{Y3},
the present author classifies a class of  higher dimensional
singular $(\alpha,\beta)$-metrics $F=\alpha\phi(\beta/\alpha)$
 which are Douglasian and locally projectively flat
respectively, where $\phi(s)$ satisfies the condition (\ref{j2}).
In this paper we will solve the singular case under  the condition
(\ref{j2}) in two-dimensional case, which shows  below that the
singular case is quite different form the regular condition
$\phi(0)=1$ (cf. \cite{Y1}).

\begin{thm}\label{Th1}
  Let $F=\alpha \phi(s)$, $s=\beta/\alpha$, be a two-dimensional
  $(\alpha,\beta)$-metric on an open subset $U\subset R^2$, where $\phi$ satisfies (\ref{j2}).
   Suppose $db\ne0$ in $U$ and that $\beta$ is not parallel with respect to
     $\alpha$.
  If $F$ is a Douglas metric, or locally projectively flat, then $F$ must be in the following
  form
  \be\label{j002}
F=c\bar{\beta}+\bar{\beta}^m\bar{\alpha}^{1-m},\ \
(\bar{\alpha}:=\sqrt{\alpha^2+k\beta^2}, \ \bar{\beta}:=\beta),
 \ee
 where $c,k$ are constant. Note that $\bar{\alpha}$ is Riemannian if
 $k>-1/b^2$.
\end{thm}

If $b=constant$ in Theorem
   \ref{Th1}, there are other classes for the metric $F$ (see Theorem \ref{th1} and Theorem \ref{th3}
   below). Theorem \ref{Th1} also holds if $n\ge 3$, but there is
   much difference between $n=2$ and $n\ge 3$ when we determine
   the local structures of $F$ in (\ref{j002}) which is Douglasian
   or locally projectively flat (cf. \cite{Y3}).

Theorem \ref{Th1}  naturally induces an important class of
singular Finsler metric---$m$-Kropina metric
$F=\beta^m\alpha^{1-m}$. When $m=-1$, $F=\alpha^2/\beta$ is called
a Kropina metric. There have been some research papers on Kropina
metrics (\cite{RR} \cite{YN} \cite{YO}). In \cite{ShY}, the
present author and Z. Shen characterize $m$-Kropina metrics which
are weakly Einsteinian.

\bigskip

Next we   determine the local structure of the
 metric $F=c\beta+\beta^m\alpha^{1-m}$ which is Douglasian and locally projectively flat respectively.
 For projectively flat case, it is hard to deal with $m=-1$ and $c\ne 0,m=-3$. The method is to apply an interesting deformation on
 $\alpha$ and $\beta$ which is defined by
  \be\label{cr71}
\widetilde{\alpha}:=b^m\alpha, \ \
\widetilde{\beta}:=b^{m-1}\beta,
 \ee
The deformation (\ref{cr71}) first appears in  \cite{ShY} for the
research on weakly Einstein $m$-Kropina metrics. It also appears
in \cite{Y3}. It is very useful for $m$-Kropina metrics.
Obviously, if $F$ is an $m$-Kropina metric, then $F$ keeps
formally unchanged, namely,
$$F=\beta^m\alpha^{1-m}=\widetilde{\beta}^m\widetilde{\alpha}^{1-m}.$$
Further, $\widetilde{\beta}$ has unit length with respect to
$\widetilde{\alpha}$, that is,
$||\widetilde{\beta}||_{\widetilde{\alpha}}=1$.

\begin{thm}\label{th001}
 Let $F=c\beta+\beta^m\alpha^{1-m}$ be a two-dimensional Douglas $(\alpha,\beta)$-metric,
 where $c,m$ are constant with $m\ne
 0, 1$. Then  we have the following  cases:

 \ben

 \item[{\rm (i)}] {\rm($m=-1$)} $\alpha$ and $\beta$ can be
 arbitrary, namely, a two-dimensional metric $F=c\beta+\alpha^2/\beta$ is always a Douglas
 metric.

 \item[{\rm (ii)}] {\rm($m=-3$)}
 $\alpha$ and $\beta$ can be locally written as
  \beq
 \alpha^2&=&\frac{B^3}{u^2+v^2}\big\{(y^1)^2+(y^2)^2\big\}-\frac{3(5+3cB^2)(1+cB^2)}{B}\beta^2,\label{yc1}\\
 \beta&=&\frac{B^2}{(4+3cB^2)(u^2+v^2)}(uy^1+vy^2),\label{yc2}
  \eeq
  where $B=B(x)>0,
   u=u(x),v=v(x)$ are  scalar functions such that
   \be\label{yc4}
    f(z)=u+iv, \ \ z=x^1+ix^2
    \ee
    is a complex analytic function.

\item[{\rm (iii)}] {\rm($m\ne -1$; $c=0$)}  $F$ can be written as
 $F=\widetilde{\alpha}^{1-m}\widetilde{\beta}^m$,
  where $\widetilde{\alpha}$ is flat and $\widetilde{\beta}$
  is parallel with respect to $\widetilde{\alpha}$, and thus $\widetilde{\alpha}$ and
  $\widetilde{\beta}$ can be locally written as
 \be\label{ycw16}
\widetilde{\alpha}=|y|,\ \ \widetilde{\beta}=y^1.
 \ee
Further $\alpha,\beta$ are related
    with $\widetilde{\alpha},\widetilde{\beta}$ by
    \be\label{ycw17}
  \alpha=\eta^{\frac{m}{m-1}}\widetilde{\alpha}, \ \ \ \beta=\eta
  \widetilde{\beta},
    \ee
where $\eta=\eta(x)>0$ is a scalar function. Obviously, $F$ is
locally Minkowskian.

\item[{\rm (iv)}] {\rm($m\ne -1$; $c\ne 0,m\ne -3$)} $F$ can be
written as $
  F=c\eta\widetilde{\beta}+\widetilde{\beta}^m\widetilde{\alpha}^{1-m},
 $
  where (\ref{ycw16}) and (\ref{ycw17}) hold with
  $\eta=\eta(x^1)>0$.
 \een
\end{thm}

In Theorem \ref{th001}, we cannot give a detailed description for
the local structures of $F$ in higher dimensions since
$\widetilde{\alpha}$ cannot be determined in this case, and
Theorem \ref{th001}(i) does not hold in higher dimensions either
(cf. \cite{Y3}). Theorem \ref{th001}(ii) and (iii) give two
representations for the local structure of $F=\alpha^4/\beta^3$.
We will prove Theorem \ref{th001}(ii)-(iv) by aid of the result in
\cite{Y} (also see \cite{Y2}). In Theorem \ref{th001}(ii), the
metric is determined by the triple parametric functions $B,u,v$,
where $u,v$ are a pair of conjugate harmonious functions.

\begin{thm}\label{th04}
 Let $F=c\beta+\beta^m\alpha^{1-m}$ be a two-dimensional
 $(\alpha,\beta)$-metric, where $c,m$ are constant with $m\ne
 0,\pm 1$ and $c=0$ if $m=-3$. Then  $F$ is Douglasian if and only if $F$ is locally
 projectively flat. In this case, $F$ is Berwaldian, or
 locally Minkowskian if and only if  $c=0$ or $\eta=constant$ in
 (\ref{ycw17}); and here $\eta=constant$ implies $\alpha$ is flat
 and $\beta$ is parallel.
\end{thm}

Theorem \ref{th04} does not hold in higher dimensions (cf.
\cite{Y3}). For $m=-1$, a two-dimensional metric in the form
$F=c\beta+\alpha^2/\beta$ is possibly Not locally projectively
flat, although it is always Douglasian. For $m=-3$, by (\ref{yc1})
and (\ref{yc2}), we can construct two-dimensional metrics in the
form $F=c\beta+\alpha^4/\beta^3$ with $c\ne 0$ which are Douglaian
but Not locally projectively flat. See the examples in the last
section. Besides, the local structure has been determined in
\cite{Y4} for $F=c\beta+\alpha^2/\beta$, or
$F=c\beta+\alpha^4/\beta^3$ which is locally projectively flat
with constant flag curvature in dimension $n\ge 2$.

\bigskip

\noindent{\bf Open Problem:} Determine the local structure of the
two-dimensional metric $F=c\beta+\alpha^2/\beta$, or
$F=c\beta+\alpha^4/\beta^3$ ($c\ne0$) which is locallly
projectively flat.

\section{Preliminaries}

Let $F=F(x,y)$ be a Finsler metric on an $n$-dimensional manifold
$M$.
      In local coordinates, the spray coefficients $G^i$ are
      defined by
 \beq \label{G1}
 G^i:=\frac{1}{4}g^{il}\big \{[F^2]_{x^ky^l}y^k-[F^2]_{x^l}\big \}.
 \eeq
If $F$ is a Douglas metric, then  $G^i$   are in the following
form:
 \beq\label{G2}
   G^i=\frac{1}{2}\Gamma_{jk}^i(x)y^jy^k+P(x,y)y^i,
 \eeq
where $\Gamma_{jk}^i(x)$ are local functions on $M$ and $P(x,y)$
is a local positively homogeneous function of degree one in $y$.
It is easy to see that $F$
 is a Douglas metric if and only if $G^iy^j-G^jy^i$ is a
 homogeneous polynomial in $(y^i)$ of degree three, which by
 (\ref{G2}) can be written as (\cite{BaMa}),
  $$G^iy^j-G^jy^i=\frac{1}{2}(\Gamma^i_{kl}y^j-\Gamma^j_{kl}y^i)y^ky^l.$$

According to G. Hamel's result, a Finsler metric $F$ is
projectively flat in $U$ if and only if
 \be\label{01}
 F_{x^my^l}y^m-F_{x^l}=0.
 \ee
 The above formula implies that $G^i=Py^i$ with $P$ given by
 \be\label{02}
  P=\frac{F_{x^m}y^m}{2F}.
  \ee

Consider an $(\alpha,\beta)$-metric $F =\alpha \phi
(\beta/\alpha)$.   The spray coefficients $G^i_{\alpha}$ of
$\alpha$
 are given by
  $$G^i_{\alpha}=\frac{1}{4}a^{il}\big \{[\alpha^2]_{x^ky^l}y^k-[\alpha^2]_{x^l}\big
  \}.$$
  Let $\nabla \beta = b_{i|j} y^i dx^j$  denote
the covariant derivatives of $\beta$ with respect to $\alpha$ and
define
 $$r_{ij}:=\frac{1}{2}(b_{i|j}+b_{j|i}),\ \ s_{ij}:=\frac{1}{2}(b_{i|j}-b_{j|i}),\ \
 r_j:=b^ir_{ij},\ \ s_j:=b^is_{ij},\ \ s^i:=a^{ik}s_k,$$
 where $b^i:=a^{ij}b_j$ and $(a^{ij})$ is the inverse of
 $(a_{ij})$.
 By (\ref{G1}) again, the spray coefficients $G^i$ of $F$
are given by:
  \be\label{y20}
  G^i=G^i_{\alpha}+\alpha Q s^i_0+\alpha^{-1}\Theta (-2\alpha Q
  s_0+r_{00})y^i+\Psi (-2\alpha Q s_0+r_{00})b^i,
  \ee
where $s^i_j=a^{ik}s_{kj}, s^i_0=s^i_ky^k,
s_i=b^ks_{ki},s_0=s_iy^i$, and
 $$
  Q:=\frac{\phi'}{\phi-s\phi'},\ \
  \Theta:=\frac{Q-sQ'}{2\Delta},\ \
  \Psi:=\frac{Q'}{2\Delta},\ \ \Delta:=1+sQ+(b^2-s^2)Q'.
 $$

By (\ref{y20}) one  can see that $F=\alpha\phi(\beta/\alpha)$ is a
Douglas metric if and only if
 \be\label{y21}
 \alpha Q (s^i_0y^j-s^j_0y^i)+\Psi (-2\alpha
 Qs_0+r_{00})(b^iy^j-b^jy^i)=\frac{1}{2}(G^i_{kl}y^j-G^j_{kl}y^i)y^ky^l,
\ee
  where $G^i_{kl}:=\Gamma^i_{kl}-\gamma^i_{kl}$,  $\Gamma^i_{kl} $ are given in (\ref{G2}) and $
  \gamma^i_{kl}:=\pa^2G^i_{\alpha}/\pa y^k\pa y^l.$

Further,  $F=\alpha\phi(\beta/\alpha)$ is projectively flat on
$U\subset R^n$ if and only if \be
(a_{ml}\alpha^2-y_my_l)G^m_{\alpha}+\alpha^3Qs_{l0}+\Psi\alpha(-2\alpha
 Qs_0+r_{00})(\alpha b_l-sy_l)=0,\label{y21*}
\ee
 where $y_l=a_{ml}y^m$.

\section{Equations in a Special Coordinate System}\label{s3}

In order to prove Theorems \ref{th1} and \ref{th3} below, one has
to simplify  (\ref{y21}) and (\ref{y21*}). The main technique is
to fix a point and choose a   special  coordinate
 system $(s,y^a)$  as  in \cite{Shen1} \cite{Shen2}.

Fix an arbitrary point $x\in M$ and take  an orthogonal basis
  $\{e_i\}$ at $x$ such that
   $$\alpha=\sqrt{\sum_{i=1}^n(y^i)^2},\ \ \beta=by^1.$$
Then we change coordinates $(y^i)$ to $(s, y^a)$ such that
  $$\alpha=\frac{b}{\sqrt{b^2-s^2}}\bar{\alpha},\ \
  \beta=\frac{bs}{\sqrt{b^2-s^2}}\bar{\alpha}, $$
where $\bar{\alpha}=\sqrt{\sum_{a=2}^n(y^a)^2}$. Let
 $$\bar{r}_{10}:=r_{1a}y^a, \ \ \bar{r}_{00}:=r_{ab}y^ay^b, \ \
 \bar{s}_0:=s_ay^a.$$
We have $\bar{s}_0=b\bar{s}_{10},s_1=bs_{11}=0$.

The following two lemmas are needed and   are trivial:

\begin{lem}\label{lem5.1}
 Under the special local coordinate system at $x$, if
 $b=constant$, then $r_{12}+s_{12}=0$ at $x$.
\end{lem}

\begin{lem}
  For $n\ge 2$, suppose $p+q\bar{\alpha}=0$, where
  $p=p(\bar{y})$ and $q=q(\bar{y})$ are homogeneous
  polynomials in $\bar{y}=(y^a)$, then $p=0,q=0$.
 \end{lem}

By the above coordinate $(s,y^a)$ and using (\ref{y21}) and
(\ref{y21*}), it follows from \cite{LSS} \cite{Shen1} we have the
following results:

\begin{prop}\label{prop3.4}
  $(n=2)$ An  $(\alpha,\beta)$-metric
 $F=\alpha\phi(\beta/\alpha)$ is a Douglas metric if and only if
 \be\label{y27}
  \frac{s^2}{2(b^2-s^2)}(G^1_{11}-G^2_{12}-G^2_{21})+\frac{1}{2}G^1_{22}=
   b\Psi (\frac{r_{11}s^2}{b^2-s^2}+r_{22}),
 \ee
 \be\label{y28}
 (-\frac{s^2}{b^2-s^2}+2\Psi
 b^2-1)bQs_{12}-2b\Psi
 r_{12}s=\frac{G^2_{11}}{2(b^2-s^2)}s^3+\frac{1}{2}(G^2_{22}-G^1_{12}-G^1_{21})s,
 \ee
 where $G^i_{jk}$ are defined in (\ref{y21}).
 \end{prop}

\begin{prop}\label{prop3.5}
  $(n=2)$ An $(\alpha,\beta)$-metric
 $F=\alpha\phi(\beta/\alpha)$ is projectively flat if and only if
 \be\label{g26}
  \frac{s^2}{2(b^2-s^2)}(-\widetilde{G}^1_{11}+2\widetilde{G}^2_{12})-\frac{1}{2}\widetilde{G}^1_{22}=
   b\Psi (\frac{s^2}{b^2-s^2}r_{11}+r_{22}),
 \ee
 \be\label{g27}
 \frac{1}{b^2-s^2}\big[2\Psi
 (b^2-s^2)-1\big]b^3Qs_{12}-2b\Psi r_{12}s=-\frac{\widetilde{G}^2_{11}}{2(b^2-s^2)}s^3+
 \frac{1}{2}(-\widetilde{G}^2_{22}+2\widetilde{G}^1_{12})s,
 \ee
 where $\widetilde{G}^i_{jk}:=\frac{\pa^2G_{\alpha}^i}{\pa y^j\pa y^k}$ are
 the connection coefficients of $\alpha$.
 \end{prop}

Comparing (\ref{y27}) and (\ref{g26}), (\ref{y28}) and
(\ref{g27}), it is easy to see that if we can solve $G^i_{jk}$
from (\ref{y27}) and (\ref{y28}), then we can solve
$\widetilde{G}^i_{jk}$ from (\ref{g26}) and (\ref{g27}). In the
following we only consider (\ref{y27}) and (\ref{y28}), from which
we will solve $G^i_{jk}$.

\section{Douglas $(\alpha,\beta)$-metrics}

In this section, we characterize a class of two-dimensional
$(\alpha,\beta)$-metrics  which are Douglasian. We have the
following theorem.

 \begin{thm}\label{th1}
  Let $F=\alpha \phi(s)$, $s=\beta/\alpha$, be a two-dimensional
  $(\alpha,\beta)$-metric on an open subset $U\subset R^2$, where $\phi$ satisfies (\ref{j2}).
   Suppose that $\beta$ is not parallel with respect to
     $\alpha$.
  Then $F$ is a Douglas metric if and only if one of the following
  cases holds:
  \ben
 \item[{\rm (i)}]  $\phi$ is given by
 \be\label{ygjcw}
\phi(s)=cs+\frac{1}{s},
 \ee
and $\alpha$ and $\beta$ are
 arbitrary, where $c$ is a constant.

  \item[{\rm (ii)}]  $\phi$ and $\beta$ satisfy
\beq
    \phi(s)&=&k_1s+\frac{2k_2}{s}+\frac{1}{s^3},\label{yg5} \\
   r_{ij}&=&-2\tau \big\{3b^2a_{ij}+(k_2b^2-2)b_ib_j\big\}+\frac{(3k_1+k_2^2)b^4-4}{8b^2(1+k_2b^2)}(b_is_j+b_js_i),\label{yg6}
   \eeq
   where $\tau=\tau(x)$ is a scalar function and $k_1,k_2$ are
constant satisfying $1+k_2b^2\ne 0$.

\item[{\rm (iii)}]
  $\phi$ and $\beta$ satisfy
   \beq
    \phi(s)&=&k_1s+s^m(1+k_2 s^2)^{\frac{1-m}{2}},\label{y5} \\
     b_{i|j}&=&2\tau \big\{mb^2a_{ij}-(m+1+k_2b^2)b_ib_j\big\}, \label{y6}
   \eeq
   where $\tau=\tau(x)$ is a scalar function and $k_1,k_2$ are
constant.

  \item[{\rm (iv)}]  $\phi$ and $\beta$ satisfy
   \beq
    \phi(s)&=&s^m(1+k s^2)^{\frac{1-m}{2}},\label{y16}\\
   r_{ij}&=&2\tau \big\{mb^2a_{ij}-(m+1+kb^2)b_ib_j\big\}-\frac{m+1+2kb^2}{(m-1)b^2}(b_is_j+b_js_i),\label{y17}
   \eeq
  where $k$ is constant and $\tau=\tau(x)$ is a scalar.

 \item[{\rm (v)}]  $\phi$ and $\beta$ satisfy
   \beq
    \phi(s)&=&mb^2\sqrt{b^2-s^2}\int_0^s\frac{1}{(b^2-t^2)^{3/2}}
    \Big(\frac{t}{\sqrt{1-k t^2}}\Big)^{m-1}dt,\label{j4}\\
   r_{ij}&=&-\frac{1}{b^2}(b_is_j+b_js_i),\label{j5}
   \eeq
   where $k$ is a constant.
  \een
 \end{thm}

In Theorem \ref{th1} (iv), if $b=constant$, then $k=-1/b^2$ in
(\ref{y16})--(\ref{y17}), and we get
   \beq
    \phi(s)&=&s^m\big\{1-(\frac{s}{b})^2\big\}^{\frac{1-m}{2}},\label{y18}\\
    r_{ij}&=&2\bar{\tau}(b^2a_{ij}-b_ib_j)-\frac{1}{b^2}(b_is_j+b_js_i),\label{y19}
   \eeq
where $\bar{\tau}:=m\tau$. Note that if $n=2$, (\ref{y19}) is
equivalent to $b=constant$ (see \cite{LS2}).  When $k_1=k_2^2$,
Theorem \ref{th1} (ii) is a special case of Theorem \ref{th1}
(iv).

\subsection{ $(r_{11},r_{22})\ne (0,0)$}

 {\bf Step I: } We first consider the equation (\ref{y27}). Since $(r_{11},r_{22})\ne (0,0)$,
 (\ref{y27}) can be written as
 \be
  2\Psi=\frac{\lambda s^2+\mu (b^2-s^2)}{\delta s^2+\eta
  (b^2-s^2)},\label{c2}
 \ee
 where
 $\lambda=\lambda(x),\mu=\mu(x),\delta=\delta(x),\eta=\eta(x)$ are
 scalar functions.  Since $\phi$ satisfies (\ref{j2}), $F$ is not of Randers
 type and we have $\lambda\eta-\mu\delta\ne 0$. Then by (\ref{y27}) and (\ref{c2}), for some scalar $\bar{\tau}=\bar{\tau}(x)$,
there hold (see also \cite{LSS})
 \be\label{y30}
 r_{11}=2b^2\delta \bar{\tau},\ \ r_{22}=2b^2\eta \bar{\tau},
 \ee
 \be\label{y0030}
 G^1_{11}=G^2_{12}+G^2_{21}+2\lambda b^3\bar{\tau},\ \ G^1_{22}=2\mu
 b^3\bar{\tau}.
 \ee
 Rewrite (\ref{c2}) as
follows
 \be\label{y56}
 [\delta s^2+\eta (b^2-s^2)]\phi''=[\lambda s^2+\mu
 (b^2-s^2)][\phi-s\phi'+(b^2-s^2)\phi''].
 \ee
Plug
 \be\label{y0056}
 \phi(s)=cs+s^m(1+a_{m+1}s+a_{m+2}s^2+a_{m+3}s^3+a_{m+4}s^4+a_{m+5}s^5)+o(s^{m+5})
 \ee
into (\ref{y56}). Let $p_i$ be the coefficients of $s^i$ in
(\ref{y56}). First $p_{m-2}=0$ gives
 \be\label{y57}
 \eta=\mu b^2.
 \ee
  Plugging (\ref{y57}) into  $p_m=0$
yields
 \be\label{y58}
 \delta=\lambda b^2-\frac{m+1}{m}\mu b^2.
 \ee

\noindent{\bf Case A. } Assume $m=-1$. Plug (\ref{y57}),
(\ref{y58}) and $m=-1$ into (\ref{y56}) and then we get
 $$s^2\phi''+s\phi'-\phi=0,$$
 whose solution is given by (\ref{ygjcw}). By (\ref{y30}),
(\ref{y0030}), (\ref{y57}) and (\ref{y58}) we obtain
 \be\label{cr1}
G^1_{11}=G^2_{12}+G^2_{21}+\frac{r_{11}}{b},\ \
G^1_{22}=\frac{r_{22}}{b}.
 \ee

\noindent{\bf Case B. } Assume $m\ne -1$. Plugging (\ref{y57}) and
(\ref{y58}) into  $p_{m+2}=0$ yields
 \be\label{y59}
  \lambda=[m(m-1)+2a_{m+2}b^2]\epsilon, \ \ \mu=m(m-1)\epsilon,
 \ee
where $\epsilon=\epsilon(x)\ne 0$ is a scalar. It is easy to see
that
 \be\label{y059}
 \lambda \eta-\mu\delta=m(m+1)(m-1)^2b^2\epsilon^2\ne 0.
 \ee
 Plug (\ref{y57}), (\ref{y58}) and (\ref{y59}) into
(\ref{c2}) and we get
 \be\label{y60}
 2\Psi=\frac{\phi''}{\phi-s\phi'+(b^2-s^2)\phi''}=\frac{m(m-1)+2a_{m+2}s^2}{m(m-1)b^2+(1-m^2+2a_{m+2}b^2)s^2}.
 \ee

\bigskip

\noindent{\bf Step II: } Next we solve the equation (\ref{y28}).
Put $\xi:=G^1_{12}+G^1_{21}-G^2_{22}$.\

\bigskip

 \noindent{\bf Case A. } Assume $m=-1$. Plug (\ref{ygjcw}) into
 (\ref{y28}) and we have
  \be\label{ycw55}
 \xi=\frac{2r_{12}-s_{12}}{b},\ \ G^2_{11}=\frac{1-cb^2}{b}s_{12}.
  \ee
Now we have seen that a two dimensional metric
$F=k\beta+\alpha^2/\beta$ is always  Douglasian.

\bigskip

 \noindent{\bf Case B. } Assume $m\ne -1$. Plug (\ref{y60}) into (\ref{y28}) and we have
 \beq
 &&0=-2b(b^2-s^2)[m(m-1)+2a_{m+2}s^2]r_{12}(\phi-s\phi')-2b^3\phi's(1-m+2a_{m+2}s^2)s_{12}\nonumber\\
 &&\ \ \ \ \ +[(1-m^2+2a_{m+2}b^2)s^2+m(m-1)b^2][(b^2-s^2)\xi-G^2_{11}s^2](\phi-s\phi').\label{y62}
 \eeq
Plug the expansion as in (\ref{y0056}) into (\ref{y62}). Let $p_i$
denote the coefficient of $s^i$ in (\ref{y62}).
 For convenience, we put
  \be\label{y066}
  a_{m+2}=\frac{1}{2}c_1,\ \
  a_{m+4}=\frac{1}{8}\frac{(m^2+5m+4)c_1^2-2c_2}{m(m-1)}.
  \ee
Note that in the next computation, when $c\ne 0$ and $m=3$, $p_m$
has different result from that for $m\ne 3$, but there is no
effect on the final result. So we only consider $m\ne 3$ in the
computation for $p_m$. Solving the system $p_m=0,p_{m+2}=0$ and
$p_{m+4}=0$ yields the following two cases:
 \ben
  \item[{\rm (i)}] If
  \be\label{y070}
   m-1-c_1b^2\ne 0,
  \ee
then we have
 \be\label{y67}
 r_{12}=\frac{c_2b^4-(m^2-1)(m+3)c_1b^2+(m^2-1)^2}{(m+1)(m-1)^2(1-m+c_1b^2)}s_{12},
 \ee
 \be\label{y68}
  G^2_{11}=2\frac{\big\{(m+2)c_1^2-c_2\big\}b^4+m(m^2-1)c_1b^2-m^2(m-1)^2}
  {m(m-1)^2b(1-m+c_1b^2)}s_{12},
 \ee
 \be\label{y69}
 \xi=2\frac{c_2b^4-(m^2-1)(m+2)c_1b^2+m(m+1)(m-1)^2}{(m+1)(m-1)^2b(1-m+c_1b^2)}s_{12}.
 \ee

 \item[{\rm (ii)}] If
 \be
   m-1-c_1b^2=0,\label{y70}
  \ee
 then we get
  \ben
    \item[{\rm (iia)}] If $s_{12}=0$, then $r_{12}=0$ by Lemma
    \ref{lem5.1} since $b=constant$.

   \item[{\rm (iib)}] If $s_{12}\ne 0$, then
   \be\label{y74}
   c_1=\frac{m-1}{b^2},\ \ c_2=2(m+1)c_1^2.
   \ee
  \een
 \een

Now we can determine $Q$ under $s_{12}\ne 0$, and $\phi$ under two
cases: $s_{12}=0$ and $s_{12}\ne 0$.

\bigskip

 \noindent {\bf Case B1:} Assume $s_{12}=0$.

 If (\ref{y070}) holds, then $r_{12}=0$ by
 (\ref{y67}). If (\ref{y70}) holds, we also get $r_{12}=0$. Then by
 (\ref{y30}), (\ref{y57}),  (\ref{y58}) and (\ref{y59}) we get the expression of
 $b_{i|j}$ given by (\ref{y6}), where we put

 \be\label{y060}
 \tau:=(m-1)\epsilon b^2\bar{\tau}.
 \ee
   By (\ref{y60}) we get
 \be\label{y61}
 \phi''=\frac{-m+k_2s^2}{(1+k_2s^2)s^2}(\phi-s\phi'),
 \ee
 where we put
 $$k_1=a_1,\ \ k_2=-2a_{m+2}/(m-1).$$
Solving the differential equation (\ref{y61}) gives (\ref{y5}).
 This class belongs to Theorem \ref{th1}(iii).

 \bigskip

\noindent {\bf Case B2:}  Assume $s_{12}\ne 0$.

 {\bf  B2(i):} Suppose (\ref{y070}). We plug (\ref{y67}),
   (\ref{y68}) and (\ref{y69}) into
  (\ref{y62}), and then we obtain
  \be\label{y75}
 Q= \frac{[(m+2)c_1^2-c_2]s^4+m(m^2-1)c_1s^2-m^2(m-1)^2}
  {m(m-1)^2s(1-m+c_1s^2)},
  \ee
By (\ref{y75}) we can get $\phi'$ and by differentiating it we get
$\phi''$. Then plugging $\phi'$ and $\phi''$ into (\ref{y60}) we
have the following two cases:

\ \ \  {\bf  B2(i)(1):} $m=-3$. In this case, (\ref{y75}) implies
(\ref{y60}). Then solving (\ref{y75}) gives (\ref{yg5}), where we
define $k_1,k_2$ by
 \be\label{cr2}
 k_1:=-\frac{c_2+c_1^2}{48},\ \ k_2:=\frac{c_1}{4}.
 \ee
By (\ref{y30}), (\ref{y57}), (\ref{y58}),
 (\ref{y59}) and (\ref{y67}) we obtain (\ref{yg6}), where $\tau=\tau(x)$ is defined by (\ref{y060}) with $m=-3$.
 Then we get Theorem \ref{th1}(ii).

\ \ \  {\bf B2(i)(2):} $m\ne -3$. In this case, we have
 \be\label{y78}
 c_2=2(m+1)c_1^2.
\ee
  Plug (\ref{y78}) into
 (\ref{y67}) and then we have
 \be\label{y79}
  r_{12}=-\frac{m+1+2kb^2}{m-1}s_{12},
 \ee
 where $k=-c_1/(m-1)$. By (\ref{y30}), (\ref{y57}), (\ref{y58}) ,
 (\ref{y59}) and (\ref{y79}) we obtain (\ref{y17}) with $\tau:=(m-1)b^2\epsilon \bar{\tau}$. Plugging
 (\ref{y78}) and $c_1=(1-m)k$ into (\ref{y75}) yields
 $$\frac{\phi'}{\phi-s\phi'}=-\frac{m+ks^2}{(m-1)s}.$$
Thus we easily get $\phi$ given by (\ref{y16}). This class belongs
to Theorem \ref{th1}(iv).

{\bf  B2(ii):} Suppose (\ref{y70}). Then $r_{12}=-s_{12}$. Plug
(\ref{y74}) into (\ref{y57}), (\ref{y58}) and (\ref{y59}), then we
get $\delta,\eta$. Plug $\delta,\eta$ into (\ref{y30}), then we
get $r_{11},r_{22}$. Plus $r_{12}=-s_{12}$ we obtain (\ref{y19})
for some scalar $\bar{\tau}=\bar{\tau}(x)$.  Similarly as above we
get $\phi$ given by (\ref{y18}). This class belongs to Theorem
\ref{th1}(iv).

   \subsection{$(r_{11},r_{22})=(0,0)$}

Since $\beta$ is not parallel and $(r_{11},r_{22})= (0,0)$, we
will see that $s_{12}\ne 0$ from
 the following proof for different cases. It follows
from   (\ref{y27}) that
\[ G^1_{22}=0, \ \ \ \ \ \ G^1_{11}=G^2_{12}+G^2_{21}.\]
Plugging the expressions of $Q$ and $\Psi$ into (\ref{y28}) yields
 \beq
 &&s(b^2-s^2)\big[2br_{12}+G^2_{11}s^2-(b^2-s^2)\xi\big]\phi''\nonumber\\
 &&\ \ \ \ +
 s\big[G^2_{11}s^2-(b^2-s^2)\xi\big](\phi-s\phi')+2b^3s_{12}\phi'=0,\label{yy49}
 \eeq
where $\xi:=G^1_{12}+G^1_{21}-G^2_{22}$.

Plug
 $$\phi=a_1s+s^m\big(1+a_{m+1}s+a_{m+2}s^2+a_{m+3}s^3+a_{m+4}s^4+\cdots\big)$$
into (\ref{yy49}) and let $p_i$ be the coefficient of $s^i$ in
(\ref{yy49}). All $p_i$'s are zero.

 By $p_{m-1}=0$ we have
 \be\label{w61}
 r_{12}=\frac{1}{2}b\xi-\frac{1}{m-1}s_{12}.
 \ee
Plugging (\ref{w61}) into $p_{m+1}=0$ yields
 \be\label{w62}
 G^2_{11}=-\frac{m+1}{m}\xi+\frac{2}{m-1}\Big\{\frac{2(m+2)ba_{m+2}}{m(m-1)}-\frac{1}{b}\Big\}s_{12}.
 \ee

\noindent {\bf Case I:} Assume $b\ne constant$. We will get
Theorem \ref{th1}(i), (ii) and (iv) in a special case.

\bigskip

\noindent {\bf Case IA:} $m=-1$. By the discussion of the two
cases (1) and (2) below we get $\phi(s)$ given by (\ref{ygjcw})
for some constant $c$.

{\bf (1).} Assume $a_{2k}\ne 0$ for some minimal integer $k\ge 0$.
Plugging $a_{2k-2}=0$ into $p_{2k+1}=0$ gives
 \be\label{w062}
 \xi=\Big\{\frac{2ka_1}{1+2k}-\frac{2(1+k)(3+2k)a_{2k+2}}{(2k+1)(2k-1)a_{2k}}\Big\}bs_{12}.
 \ee
Then substitute (\ref{w61}), (\ref{w62}) and (\ref{w062}) into
(\ref{yy49}), and using $b\ne constant$ we obtain $\phi(s)$ given
by (\ref{ygjcw}) for some constant $c$. Thus  $a_{2k}=0$ for all
integers $k\ge 0$ by (\ref{ygjcw}). This contradicts with our
assumption.

{\bf (2).} Assume $a_{2k+1}\ne 0$ for some minimal integer $k\ge
1$. If $k=1$, we get $s_{12}=0$ by $p_2=0$. If $k=2$, we get
$s_{12}=0$ by plugging $a_3=0$ into $p_4=0$. If $k\ge 3$, we get
$s_{12}=0$ by plugging $a_{2k-3}=0$ and $a_{2k-1}=0$ into
$p_{2k}=0$. Now substitute (\ref{w61}), (\ref{w62}) and $s_{12}=0$
into (\ref{yy49}), and then we get
 \be\label{w0062}
 \xi s(b^2-s^2)(s^2\phi''+s\phi'-\phi)=0.
\ee
 If $\xi=0$, then we have $r_{12}=0$ by (\ref{w61}) and
 $s_{12}=0$. Thus $\xi\ne 0$. Then by (\ref{w0062})
 we obtain $\phi(s)$ given by
(\ref{ygjcw}) for some constant $c$. So  $a_{2k+1}=0$ for all
integers $k\ge 1$ by (\ref{ygjcw}). This again contradicts with
our assumption.

\bigskip

\noindent {\bf Case IB:} $m\ne -1$. By aid of (\ref{y066}),
plugging (\ref{w61}) and (\ref{w62}) into $p_{m+3}=0$ yields
 \be\label{w64}
 \xi=\frac{2\Big\{\big[(m+4)c_2-2(m+1)(m+3)c_1^2\big]b^4-(m^2-1)\big[(m+2)c_1b^2-m(m-1)\big]\Big\}}{(m+1)(m-1)^2b(1-m+c_1b^2)}s_{12}.
 \ee
By (\ref{w61}) and (\ref{w64}) we conclude that if $s_{12}=0$,
then we have $r_{12}=0$. So in this case, we have $s_{12}\ne 0$.
Plug (\ref{w61}), (\ref{w62}) and (\ref{w64}) into (\ref{yy49})
and then we obtain an equation in the form
 $$f_0+f_2b^2+f_4b^4=0,$$
where $f_0,f_2$ and $f_4$ are ODEs about $\phi$. Since $b\ne
constant$, we have $f_0=0,f_2=0,f_4=0$. This system is equivalent
to $f_0=0,f_2=0$.

{\bf (1).} $m=-3$. In this case, solving the system $f_0=0,f_2=0$
gives (\ref{yg5}) for some constants $k_1,k_2$. Meanwhile we get
(\ref{yg6}) with $\tau=0$ by (\ref{w61}), (\ref{w64}), $r_{11}=0$
and $r_{22}=0$.

{\bf (2).} $m=-4$. We get $\phi$ given by (\ref{y16}) with $m=-4$.

{\bf (3).} $m\ne-3,-4$. Solving the system $f_0=0,f_2=0$ we can
first show that
 \be\label{w65}
 c_2=2(m+1)c_1^2.
 \ee
Plugging (\ref{w65}) into the system $f_0=0,f_2=0$ again we get
the solution of $\phi$ given by (\ref{y16}) with $k=-c_1/(m-1)$.

If $m\ne -3$, plug (\ref{w65}) into (\ref{w61}), (\ref{w62}) and
(\ref{w64}) and we have
 \be\label{w66}
 r_{12}=\frac{1+m+2kb^2}{1-m}s_{12},\ \
 G^2_{11}=\frac{2(m+kb^2)}{(m-1)b}s_{12},\ \
 \xi=\frac{2(m+2kb^2)}{(1-m)b}s_{12}.
 \ee
Now by the expression of $r_{12}$ in (\ref{w66}) we get
(\ref{y17}) with $\tau=0$.

 \bigskip

\noindent {\bf Case II:} Assume $b= constant$. We will show this
case gives the class Theorem \ref{th1}(v).

Since $r_{12}+s_{12}=0$, it follows from (\ref{w61}) that
 \be\label{w67}
 \xi=\frac{2(m-2)}{(1-m)b}s_{12}.
 \ee

If $m=-2$, substituting (\ref{w67}) and (\ref{w62}) into
(\ref{yy49}) yields (\ref{w69}) with $m=-2$ and $k=1/b^2$.

 If $m\ne -2$, plug (\ref{w67}) into (\ref{w62}) and we have
  \be\label{w68}
 G^2_{11}=\frac{2(m-2+kb^2)}{(m-1)b}s_{12},
  \ee
 where we put
  $$a_{m+2}=\frac{(m-1)(2+mkb^2)}{2(m+2)b^2}.$$
Now plugging $r_{12}=-s_{12}$, (\ref{w67}) and (\ref{w68}) into
(\ref{yy49}) yields
 \be\label{w69}
 \frac{\phi-s\phi'+(b^2-s^2)\phi''}{s\phi+(b^2-s^2)\phi'}=\frac{m-1}{s(1-ks^2)}.
 \ee
 Let
\[ \Phi:= s \phi(s)+(b^2-s^2)\phi'(s).\]
Then (\ref{w69}) becomes
\[ \frac{\Phi'}{\Phi} = \frac{m-1}{s(1-ks^2)}.\]
We get
\[ \Phi = c \Big ( \frac{s}{\sqrt{1-ks^2}} \Big )^{m-1},\]
where $c$ is a constant. Then we can easily get
 $\phi$ given by (\ref{j4}). And (\ref{j5})
follows from $r_{11}=0,r_{22}=0$ and $r_{12}=-s_{12}$.

\section{Projectively flat $(\alpha,\beta)$-metrics}
In this section, we characterize a two-dimensional
$(\alpha,\beta)$-metric $F=\alpha\phi(\beta/\alpha)$ satisfying
(\ref{j2}) which is projectively flat. We have the following
theorem.

\begin{thm}\label{th3}
  Let $F=\alpha \phi(s)$, $s=\beta/\alpha$, be an
  $(\alpha,\beta)$-metric on an open subset $U\subset R^2$.
   Suppose that $\beta$ is not parallel with respect to
     $\alpha$ and $\phi$ satisfies (\ref{j2}). Let $G^i_{\alpha}$  be the spray coefficients of $\alpha$.
  Then $F$ is  projectively flat in $U$ with
 $G^i=P(x,y)y^i$ if and only if one of the following cases holds:

 \ben
 \item[{\rm (i)}] $\phi(s)$  satisfy
  (\ref{ygjcw}), and  $G^i_{\alpha}$
  satisfy
  \be\label{w001}
G^i_{\alpha}=\rho
y^i-\frac{r_{00}}{2b^2}b^i-\frac{\alpha^2-c\beta^2}{2b^2}s^i.
 \ee
In this case, the projective factor $P$ is given by
 \be\label{w0001}
 P=\rho-\frac{1}{b^2(\alpha^2+c\beta^2)}\Big\{(\alpha^2-c\beta^2)s_0+r_{00}\beta\Big\}.
 \ee

\item[{\rm (ii)}]  $\phi(s)$ and $\beta$ satisfy (\ref{yg5}) and
  (\ref{yg6}), and  $G^i_{\alpha}$ satisfy
  \be\label{cw002}
G^i_{\alpha}=\rho
y^i+\tau(3\alpha^2+k_2\beta^2)b^i+\Big\{\frac{k_1-k_2^2}{8(1+k_2b^2)}(3b^2\alpha^2-\beta^2)
+(\frac{k_2}{2}-\frac{3}{4b^2})\alpha^2-\frac{k_2}{b^2}\beta^2\Big\}s^i.
 \ee
In this case, the projective factor $P$ is given by
 \beq\label{cw0002}
 P&=&\rho+2\tau\beta\Big\{3-\frac{2c\beta^4}{\alpha^4+c\beta^4+k_2\beta^2(2\alpha^2+k_2\beta^2)}\Big\}
 +(\frac{k_2b^2-3}{2b^2}+T)s_0,\\
 T:&=&c\frac{4\beta^2(2\beta^2-b^2\alpha^2)+3b^4(\alpha^4+c\beta^4)+k_2b^2\beta^2(6b^2\alpha^2+4\beta^2+3k_2b^2\beta^2)}
 {8b^2(1+k_2b^2)\big[\alpha^4+c\beta^4+k_2\beta^2(2\alpha^2+k_2\beta^2)\big]},\nonumber\\
 c:&=&k_1-k_2^2.\nonumber
 \eeq

  \item[{\rm (iii)}]  $\phi(s)$ and $\beta$ satisfy (\ref{y5}) and
  (\ref{y6}), and  $G^i_{\alpha}$
  satisfy
 \be\label{w1}
 G^i_{\alpha} =\rho y^i-\tau(m\alpha^2-k_2\beta^2)b^i.
  \ee
 In this case, the projective factor $P$ is given by
 \be\label{w01}
  P=\rho +\tau \alpha
  \Big\{s(-m+k_2s^2)-s^2(1+k_2s^2)\frac{\phi'}{\phi}\Big\}.
  \ee

\item[{\rm (iv)}] $\phi(s)$ and $\beta$ satisfy (\ref{y16}) and
  (\ref{y17}), and  $G^i_{\alpha}$
  satisfy
   \be\label{w3}
   G^i_{\alpha}=\rho
   y^i-\tau(m\alpha^2-k\beta^2)b^i+\frac{1}{1-m}\Big\{\big(2k+\frac{m}{b^2}\big)\alpha^2-\frac{k}{b^2}\beta^2\Big\}s^i.
   \ee
 In this case, the projective factor $P$ is given by
  \be\label{w03}
  P=\rho-2m\tau\beta-\frac{2(m+kb^2)}{(m-1)b^2}s_0.
  \ee

\item[{\rm (v)}]  $\phi(s)$ and $\beta$ satisfy (\ref{j4}) and
  (\ref{j5}), and $G^i_{\alpha}$
  satisfy
   \be\label{w4}
   G^i_{\alpha}=\rho
   y^i-\frac{(m-2)\alpha^2+k\beta^2}{(m-1)b^2}s^i.
   \ee
 In this case, the projective factor $P$ is given by
  \be\label{w04}
  P=\rho+\frac{1}{(m-1)b^2}\Big\{s(ks^2-1)\frac{\phi'}{\phi}-ks^2-m+2\Big\}s_0.
  \ee
  \een
The above function $\rho= \rho_i(x)y^i$  is a 1-form.

\end{thm}

To complete the proof of Theorem \ref{th3}, we only need to solve
$\widetilde{G}^i_{jk}$ from (\ref{g26}) and (\ref{g27}), and all
projective factors for every class in Theorem \ref{th1} when $F$
is projectively flat.

\begin{rem}
 In Theorem \ref{th3}, when $\beta$ is not closed, since $n=2$, we can express $G^i$ and
 $P$ for every class in different forms with different choices of $\rho$ using $s_{12}$ in the
 following proof. However, we can verify conversely that these different forms are equivalent to
 one another using the dimension $n=2$.
\end{rem}

\subsection{The Spray Coefficients of $\alpha$}

In this subsection  we will show the expressions of the spray
coefficients $G^i_{\alpha}$ for each class in Theorem \ref{th1}
when $F$ is projectively flat. Note that by
$\widetilde{G}^i_{jk}=\frac{\pa^2G_{\alpha}^i}{\pa y^j\pa y^k}$,
the spray $G^i_{\alpha}$ of $\alpha$ can be expressed as
 $$G^i_{\alpha}=\frac{1}{2}\widetilde{G}^i_{jk}y^jy^k.$$

\bigskip

\noindent
 {\bf Case I:}
Suppose $(r_{11},r_{22})\ne (0,0)$.

 {\bf (1).} Assume $m=-1$. It follows from (\ref{cr1}) and
 (\ref{ycw55}) that $G^i_{\alpha}$ are given by (\ref{w001}),
 where $\rho$ is defined by
  $$
 \rho:=\widetilde{G}^2_{12}y^1+(\widetilde{G}^1_{12}+\frac{r_{12}}{b})y^2.
  $$

{\bf (2).} Assume $m\ne -1$. By (\ref{y0030}) we have
 \be\label{y96}
  \widetilde{G}^1_{11}=2t_1-2\lambda b^3\bar{\tau},\ \
  \widetilde{G}^1_{22}=-2\mu
  b^3\bar{\tau}, \ \ \widetilde{G}^1_{12}=t_2, \ \
  \widetilde{G}^2_{12}=t_1.
 \ee
  By (\ref{y68}) and (\ref{y69}) we
get
 \be\label{y97}
  \widetilde{G}^2_{11}=-2\frac{\big\{(m+2)c_1^2-c_2\big\}b^4+m(m^2-1)c_1b^2-m^2(m-1)^2}
  {m(m-1)^2b(1-m+c_1b^2)}s_{12},
 \ee
 \be\label{y98}
 \widetilde{G}^2_{22}=2\widetilde{G}^1_{12}+2\frac{c_2b^4-(m^2-1)(m+2)c_1b^2+m(m+1)(m-1)^2}{(m+1)(m-1)^2b(1-m+c_1b^2)}s_{12}.
 \ee
Let $\rho=t_iy^i$ and $\tau$ be given by (\ref{y060}).
 If $\beta$ is closed ($s_{12}=0$), then it follows from
 (\ref{y57}), (\ref{y58}) and (\ref{y59}) and
 (\ref{y96})--(\ref{y98}) that (\ref{w1}) holds.
If $\beta$ is not closed, then  if $m=-3$ we get (\ref{cw002})
from (\ref{y96})--(\ref{y98}), where $k_1,k_2$ are defined by
(\ref{cr2}); if $m\ne -3$ we get (\ref{w3}) from (\ref{y78}),
(\ref{y96})--(\ref{y98}), where $k=c_1/(1-m)$,

\bigskip

 \noindent{\bf Case II:}
Suppose $(r_{11},r_{22})= (0,0)$. Then by (\ref{g26}) we get
 \be\label{g99}
  \widetilde{G}^1_{22}=0,\ \
  \widetilde{G}^1_{11}=2\widetilde{G}^2_{12}=2t_1, \ \
  \widetilde{G}^2_{12}=t_1,\ \ \widetilde{G}^1_{12}=t_2.
 \ee

{\bf (1).} $b\ne constant$. If $m=-1$, (\ref{w001}) has been
proved. If $m=-3$, we get (\ref{cw002}) with $\tau=0$. If $m\ne
-1,-3$, then it follows from (\ref{w66}) that
 \be\label{w73}
 \widetilde{G}^2_{11}=-\frac{2(m+kb^2)}{(m-1)b}s_{12},\ \
 \widetilde{G}^2_{22}=2\widetilde{G}^1_{12}+\frac{2(m+2kb^2)}{(1-m)b}s_{12}.
 \ee
 Then by
(\ref{g99}) and (\ref{w73})
 we obtain (\ref{w3}) with $\tau=0$.

{\bf (2).}  $b= constant$.   It follows from (\ref{w67}) and
(\ref{w68}) that
 \be\label{w74}
 \widetilde{G}^2_{11}=-\frac{2(m-2+kb^2)}{(m-1)b}s_{12},\ \
 \widetilde{G}^2_{22}=2\widetilde{G}^1_{12}+\frac{2(m-2)}{(1-m)b}s_{12}.
 \ee
Then by (\ref{g99}) and (\ref{w74})
 we obtain (\ref{w4}).

 \subsection{The Projective Factors}

In this subsection,  we are going to show the projective factors
for each class in Theorem \ref{th3}.

\bigskip

\noindent{\bf (1):}  We first prove (\ref{w01}).  By (\ref{y6}) we
have
 \be\label{j61}
 r_{00}=2\tau \big\{m b^2\alpha^2-(1+m+k_2b^2)\beta^2\big\}.
 \ee
Now plug $s_{i0}=0,s_0=0$ and (\ref{w1}), (\ref{y61}) and
(\ref{j61}) into (\ref{y20}), and then we obtain (\ref{w01}).

\bigskip

\noindent{\bf (2):}  For the proofs to (\ref{w03}) and
(\ref{w04}), since $\beta$ may not be closed, it is not easy to
show the projective factors in the initial {\it local projective
coordinate system} (in such a coordinate system, geodesics are
straight lines). However, it is easy to be solved by choosing
another local projective coordinate system, and then returning to
the the initial local projective coordinate system.

Fix an arbitrary point $x_o\in U\subset R^2$. By the above idea
and a suitable affine transformation, we may assume $(U,x^i)$ is a
local projective coordinate system satisfying that
$\alpha_{x_o}=\sqrt{(y^1)^2+(y^2)^2}$ and $\beta_{x_o}=by^1$. Then
at $x_o$ we have
 $$s^1=s_1=0,\ \ s^2=s_2=bs_{12}, \ \ s_0=bs_{12}y^2, \ \ b^1=b_1=b, \ \
 b^2=b_2=0.$$

Suppose (\ref{y16}) and (\ref{y17}). Then it is easy to get
 \beqn
    r_{00}&=&2\tau \big\{mb^2\alpha^2-(1+m+kb^2)\beta^2\big\}-\frac{2(1+m+2kb^2)}{(m-1)b^2}\beta s_0,\\
     s_{0}^i&=&\frac{1}{b^2}(s_0b^i-\beta s^i),\ \ \ \  Q=\frac{m+k s^2 }{
     (1-m)s},\\
 \Theta&=&\frac{ms}{(1+m+kb^2)s^2-mb^2},\ \ \Psi=\frac{ks^2-m}{2(1+m+kb^2)s^2-2mb^2}.
 \eeqn
  Plug  (\ref{w3}) and all the above expressions into (\ref{y20}), and then at $x_o$
  we see that $G^i=Py^i$, where $P$ is given by
  \be\label{y107}
  P=\rho+2m\tau by^1-\frac{2(m+kb^2)}{(m-1)b}s_{12}y^2,
  \ee
 By using
   $$bs_{12}y^2=s_0,\ \
  by^1=\beta,,$$
 it is easy to transform  (\ref{y107}) into (\ref{w03}).
   Since $x_o$ is arbitrarily chosen, (\ref{w03}) holds everywhere.

Suppose (\ref{j4}) and
  (\ref{j5}). Then it
is easy to get (where we use (\ref{w69}) in place of (\ref{j4}))
 \beqn
 s_{0}^i=\frac{1}{b^2}(s_0b^i-\beta s^i),\ \
 r_{00}=-\frac{2}{b^2}\beta s_0,
\eeqn
$$
 \Theta =\frac{s}{2(m-1)(b^2-s^2)}\Big\{s(ks^2-1)\frac{\phi'}{\phi}+2-m-ks^2\Big\},
 $$
 $$
 \Psi=\frac{\big\{ks^4+(m-2)s^2-(m-1)b^2\big\}\phi'+s(2-m-ks^2)\phi}{2(1-m)(b^2-s^2)\big\{(b^2-s^2)\phi'+s\phi\big\}}.
$$
  Plug  (\ref{w4}) and all the above expressions into (\ref{y20}), and
  then by using the relations
   $$bs_{12}y^2=s_0,\ \ (y^1)^2+(y^2)^2=\alpha^2,\ \
  by^1=\beta,\ \ \frac{\beta}{\alpha}=s,$$
we similarly get $G^i=Py^i$, where $P$ is given by (\ref{w04}).
The
  details are omitted.

  Similarly we can get the projective factors of the other classes in
  Theorem \ref{th3}. We omit the details.

\section{Proof of Theorem \ref{th001}}\label{sec7}

By Theorem \ref{th1}, we can easily characterize a two-dimensional
metric $F=c\beta+\beta^m\alpha^{1-m}$ which is Douglasian as
follows.

\begin{cor}\label{th01}
  Let $F=c\beta+\beta^m\alpha^{1-m}$ be a two-dimensional
  $(\alpha,\beta)$-metric on an open subset $U\subset R^2$, where $c,m$ are constant with $m\ne 0, 1$.
   Then for some scalar function $\tau=\tau(x)$, we have the following cases:
    \ben
  \item[{\rm (i)}] {\rm($m=-1$)} $F$ is always a Douglas metric.

\item[{\rm (ii)}] {\rm($m=-3$)} $F$  is  a Douglas metric if and
only if $\beta$ satisfies
 \be\label{ycw01}
r_{ij}=2\tau
(-3b^2a_{ij}+2b_ib_j)+\frac{3cb^4-4}{8b^2}(b_is_j+b_js_i).
 \ee

\item[{\rm (iii)}] {\rm($c\ne 0;m\ne -1,-3$)} $F$ is a Douglas
metric if and only if $\beta$ satisfies
 \be\label{y06}
   b_{i|j}=2\tau \big\{mb^2a_{ij}-(m+1)b_ib_j\big\},
   \ee

  \item[{\rm (iv)}] {\rm($c= 0;m\ne -1$)} $F$
is a Douglas metric if and only if $\beta$ satisfies
   \be\label{y017}
   r_{ij}=2\tau \big\{mb^2a_{ij}-(m+1)b_ib_j\big\}
   -\frac{m+1}{(m-1)b^2}(b_is_j+b_js_i),
   \ee
  \een
\end{cor}

\bigskip

\noindent {\it Proof of Theorem \ref{th001} :}

\bigskip

\noindent {\bf Case I:} Assume $m=-3$. In this case,
$F=c\beta+\alpha^4/\beta^3$. Define a new Riemann metric
$\widetilde{\alpha}$ and a 1-form $\widetilde{\beta}$ by
 \be\label{yg83}
\widetilde{\alpha}:=\sqrt{\xi\alpha^2+\eta\beta^2}, \ \
\widetilde{\beta}:=\beta,
 \ee
where
 $$
 \xi:=\frac{1}{b^2(4+3cb^4)}, \ \
 \eta:=\frac{3(5+8cb^4+3c^2b^8)}{b^4(4+3cb^4)}\Big\}.
 $$
Since $F=c\beta+\alpha^4/\beta^3$ is a Douglas metric, we have
(\ref{ycw01}). Now by (\ref{yg83}) and (\ref{ycw01}), a direct
computation gives
 \be\label{yg85}
 \widetilde{r}_{ij}=-\frac{16\tau b^4}{(4+3cb^4)^2}\widetilde{a}_{ij}.
 \ee
So $\widetilde{\beta}=\beta$ is a conformal form with respect to
$\widetilde{\alpha}$. Since $\widetilde{\alpha}$ is a
two-dimensional Riemann metric, we can express
$\widetilde{\alpha}$ locally as
 \be\label{yg86}
\widetilde{\alpha}:=e^{\sigma}\sqrt{(y^1)^2+(y^2)^2},
 \ee
where $\sigma=\sigma(x)$ is a scalar function. We can obtain the
local expression of  $\widetilde{\beta}=\beta$ by (\ref{yg85}) and
(\ref{yg86}) (see \cite{Y}). Then by the result in \cite{Y}, we
have
 \be\label{yg87}
 \widetilde{\beta}=\widetilde{b}_1y^1+\widetilde{b}_2y^2=e^{2\sigma}(uy^1+vy^2),
 \ee
where $u=u(x),v=v(x)$ are a pair of scalar functions such that
 $$f(z)=u+iv, \ \ z=x^1+ix^2$$
 is a complex analytic function. Finally, we give the relation
 between
  $b^2=||\beta||^2_{\alpha}$ with the triple $(\sigma,u,v)$, which can
  be done by computing the quantity
 $||\beta||^2_{\widetilde{\alpha}}$. First, by (\ref{yg86}) and
(\ref{yg87}) we get
 \be\label{yg88}
||\beta||^2_{\widetilde{\alpha}}=e^{2\sigma}(u^2+v^2).
 \ee
On the other hand, by the definition  of $\widetilde{\alpha}$ in
(\ref{yg83}), the inverse $\widetilde{a}^{ij}$ of
$\widetilde{a}_{ij}$ is given by
 $$
\widetilde{a}^{ij}=\frac{1}{\xi}\Big(a^{ij}-\frac{\eta
b^ib^j}{\xi+\eta b^2}\Big).
 $$
Now plug $\xi$ and $\eta$ into the above, and we obtain
 \be\label{yg89}
||\beta||^2_{\widetilde{\alpha}}=\widetilde{a}^{ij}b_ib_j=\frac{b^4}{4+3cb^4}.
 \ee
Thus by (\ref{yg88}) and (\ref{yg89}) we have
 \be\label{yg089}
e^{2\sigma}=\frac{b^4}{(4+3cb^4)(u^2+v^2)}.
 \ee
 Plug (\ref{yg089}) into (\ref{yg86}), (\ref{yg87})
and (\ref{yg83}) and then we get $\alpha$ and $\beta$ given by
(\ref{yc1}) and (\ref{yc2}), where we define  $B:=b^2$.

\bigskip

\noindent {\bf Case II:} Assume $m\ne -1$ and $c=0$. To prove this
case, we first show the following lemma.

\begin{lem}\label{lem001}
Let $\alpha$ be a two-dimensional Riemann metric on a manifold
$M$. If  there is a non-zero 1-form on $M$ which is parallel with
respect $\alpha$, then $\alpha$ is flat.
\end{lem}

{\it Proof :} Let $\beta$ be parallel 1-form  with respect
$\alpha$. Express $\alpha$ and $\beta$ locally as
$$
 \alpha=e^{\sigma(x)}\sqrt{(y^1)^2+(y^2)^2},\ \
 \beta=e^{2\sigma}(u(x)y^1+v(x)y^2).
$$
Since $\beta$ is also a conformal form with respect to $\alpha$,
by the result in \cite{Y} we know that $u,v$  are a pair of
conjugate harmonious functions, or equivalently, $u,v$ satisfy
 \be\label{ycw101}
 u_1=v_2,\ \ u_2=-v_1, \ \ (u_i:=u_{x^i}, \ v_i:=v_{x^i}).
 \ee
 Put $||\beta||_{\alpha}=1$. Then we have
$$
||\beta||^2_{\alpha}=e^{2\sigma}(u^2+v^2)=1.
$$
So $\alpha$ can be written as
 $$\alpha=\sqrt{\frac{(y^1)^2+(y^2)^2}{u^2+v^2}}.$$
Now using (\ref{ycw101}), it can be shown that $\alpha$ is of zero
sectional curvature.   \qed

\bigskip

Since $c=0$, we have $F=\beta^m\alpha^{1-m}$. Define a new Riemann
metric $\widetilde{\alpha}$ and a 1-form $\widetilde{\beta}$ by
(\ref{cr71}).
 Since
$F=\beta^m\alpha^{1-m}$ is a Douglas metric ($m\ne -1$)), we have
(\ref{y017}). By (\ref{cr71}) and (\ref{y017}), a direct
computation gives $\widetilde{r}_{ij}=0$. We can also give another
simple proof. Since $F$ is a Douglas metric and $F$ keeps formally
unchanged under (\ref{cr71}), by (\ref{y017}) and using
$\widetilde{b}=1$ we have
 \be\label{ygj103}
   \widetilde{r}_{ij}=2\widetilde{\tau} \big\{m\widetilde{a}_{ij}-(m+1)\widetilde{b}_i\widetilde{b}_j\big\}
   -\frac{m+1}{m-1}(\widetilde{b}_i\widetilde{s}_j+\widetilde{b}_j\widetilde{s}_i).
   \ee
Contracting (\ref{ygj103}) by $\widetilde{b}^i$ and then by
$\widetilde{b}^j$ and using $\widetilde{r}_i+\widetilde{s}_i=0$,
it is easy to get $\widetilde{r}_{ij}=0$. In case of dimension
$n=2$, given any pair $\alpha$ and $\beta$, we always have
 \be\label{ygj104}
 s_{ij}=\frac{1}{b^2}(b_is_j-b_js_i).
 \ee
 So by $\widetilde{s}_j=0$ and (\ref{ygj104}) we have
 $\widetilde{s}_{ij}=0$, which imply that $\widetilde{\beta}$ is
 closed. Thus $\widetilde{\beta}$ is parallel with respect to
 $\widetilde{\alpha}$. Therefore, by Lemma \ref{lem001}, $\widetilde{\alpha}$
 is flat. Thus $\widetilde{\alpha}$ and $\widetilde{\beta}$ can be
 locally written in the form (\ref{ycw16}).

\bigskip

\noindent {\bf Case III:} Assume $m\ne -1,c\ne 0,m\ne -3$. Since
$F=c\beta+\beta^m\alpha^{1-m}$ is a Douglas metric, we get
(\ref{y06}).
 Under the deformation (\ref{cr71}), (\ref{y06}) becomes
 $\widetilde{b}_{i|j}=0$. So by Lemma \ref{lem001}, we
 again obtain (\ref{ycw16}). By (\ref{ycw17}) and the fact that
 $\beta=\eta\widetilde{\beta}$ is closed, we get
 $\eta=\eta(x^1)$.   \qed

\begin{rem}\label{rem61}
We can give another useful local representation corresponding to
Theorem \ref{th001}(iii) and (iv). Define
 \be\label{ycw109}
\widetilde{\alpha}:=\sqrt{\frac{(y^1)^2+(y^2)^2}{u^2+v^2}},\ \
\widetilde{\beta}:=\frac{uy^1+vy^2}{u^2+v^2},
 \ee
where   $u=u(x),v=v(x)$ satisfy (\ref{ycw101}). Then
$\widetilde{\alpha}$ is flat and $\widetilde{\beta}$ is parallel
with respect to $\widetilde{\alpha}$. Put $\alpha$ and $\beta$ as
that in (\ref{ycw17}). Then $F=\beta^m\alpha^{1-m}$ with $m\ne
0,1$ is locally Minkowskian. If $\eta,u,v$ satisfy
 \be\label{ygj110}
 u_1=v_2,\ \ u_2=-v_1,\ \ \eta_1v=\eta_2u, \ \
 (u_i:=u_{x^i},etc.),
  \ee
   then $\beta$ is closed, and $F=c\beta+\beta^m\alpha^{1-m}$ with $m\ne 0,1$ is
locally projectively flat by the proof to Theorem \ref{th04}
below. However, there is no the relation $G^i=Py^i$ in such a
coordinate system.
\end{rem}

\section{Proof of Theorem   \ref{th04}}

In this section, we will prove Theorem \ref{th04} and thus the
local structure of the two-dimensional metric
$F=c\beta+\beta^m\alpha^{1-m}$ can be determined if $F$ is locally
projectively flat with $m\ne 0,\pm 1$ and $c=0$ if $m=-3$.

\bigskip

\noindent{\it Proof of Theorem   \ref{th04} :}

 Let
$F=c\beta+\beta^m\alpha^{1-m}$ be a two-dimensional Douglas
$(\alpha,\beta)$-metric,
 where $c,m$ are constant with $m\ne
 0,\pm 1$ and $c=0$ if $m=-3$.  Then by Theorem
 \ref{th001}(iii) and (iv),
$F$ can be written as
$$F=c\eta\widetilde{\beta}+\widetilde{\beta}^m\widetilde{\alpha}^{1-m},$$
  where $\eta=\eta(x^1)$ and $\widetilde{\alpha},\widetilde{\beta}$ are given by
  (\ref{ycw16}). Now we can easily verify that (\ref{01}) holds.
  So $F$ is projectively flat with $G^i=Py^i$. Further, by
  (\ref{02}) we can get the projective factor $P$ given by
 \be\label{ycw107}
 P=\frac{c\eta_1}{2F}(y^1)^2, \ \ \eta_1:=\eta_{x^1}.
 \ee
Besides,  its scalar flag curvature $K$ is given by
 \be\label{ycw108}
K=
\frac{c(y^1)^3}{2F^3}\Big\{\frac{3c\eta_1^2y^1}{2F}-\eta_{11}\Big\},\
 \ \eta_{11}:=\eta_{x^1x^1}.
\ee
 Then by (\ref{ycw107}) and (\ref{ycw108}), $F$ is Berwaldian, or
 locally Minkowskian if and only if  $c=0$ or $\eta=constant$.
 \qed

\bigskip

Note that the method applied in Theorem \ref{th001} cannot be used
to determine the local structure of the metric
$F=c\beta+\alpha^2/\beta$, or $F=c\beta+\alpha^4/\beta^3$
 ($c\ne0$) when $F$ is locally projectively flat. In this case, we
 can only obtain a general characterization by Theorem \ref{th3}, as shown in the following
 corollary.

\begin{cor}\label{th03}
  Let $F$ be a two-dimensional
  $(\alpha,\beta)$-metric. If $F=c\beta+\alpha^2/\beta$, then $F$ is locally projectively flat
  if and only if the spray $G^i_{\alpha}$ of $\alpha$ satisfy
 $$
 G^i_{\alpha}=\rho
y^i-\frac{r_{00}}{2b^2}b^i-\frac{\alpha^2}{2b^2}s^i.
 $$
 If $F=c\beta+\alpha^4/\beta^3$, then $F$ is locally projectively flat
  if and only if $\beta$ satisfies (\ref{ycw01}) and  the spray $G^i_{\alpha}$ of $\alpha$ satisfy
  $$
 G^i_{\alpha}=\rho
y^i+3\tau\alpha^2b^i+\Big\{\frac{c}{8}(3b^2\alpha^2-\beta^2)
-\frac{3\alpha^2}{4b^2}\Big\}s^i.
 $$
\end{cor}

\section{Examples}\label{sec8}

 In this section, we will construct  some examples which are
 Douglasian or projectively flat. Further, we show for the metric
 $F=c\beta+\alpha^2/\beta$, or $F=c\beta+\alpha^4/\beta^3$
 ($c\ne0$), there are examples which are Douglasian but not
 locally projectively flat.

 \begin{ex} In Remark \ref{rem61}, put
  $$u:=x^1,\ \ v:=x^2, \ \ \eta:=|x|^{1-m}.$$
  It is easy to see that $u,v,\eta$ satisfy (\ref{ygj110}), and $\alpha$ and
  $\beta$ determined by (\ref{ycw17}) and (\ref{ycw109}) are given
  by
 $$\alpha:=\frac{|y|}{|x|^{m+1}},\ \
 \beta:=\frac{\langle x,y\rangle}{|x|^{m+1}}.$$
 Then the $(\alpha,\beta)$-metric
 $F=c\beta+\beta^m\alpha^{1-m}$ is locally projectively flat, where $c,m$
 are constant with $m\ne 0,1$. But we do not have
 $G^i=Py^i$ in the present coordinate system.
 \end{ex}

 \begin{ex} In Remark \ref{rem61}, put
  $$u:=x^2,\ \ v:=-x^1, \ \ \eta:=|x|^{1-m}.$$
  It is easy to verify that $u,v,\eta$ does not satisfy the third equation in (\ref{ygj110}), and $\alpha$ and
  $\beta$ determined by (\ref{ycw17}) and (\ref{ycw109}) are given
  by
 $$\alpha:=\frac{1}{|x|^{m+1}}|y|,\ \
 \beta:=\frac{1}{|x|^{m+1}}(x^2y^1-x^1y^2).$$
Then the $(\alpha,\beta)$-metric
 $F=\beta^m\alpha^{1-m}$ is locally Minkowskian, where $c,m$
 are constant with $m\ne 0,1$. Obviously, $\beta$ is not closed.
 \end{ex}

\bigskip

Now we consider the metrics $F=c\beta+\alpha^2/\beta$, and
$F=c\beta+\alpha^4/\beta^3$. To verify the following two examples,
we need to mention the so-called $K$-curvature (\cite{Ma1}). For
an $n$-dimensional Finsler metric $F$, let $R^i_{\ k}$ be the
Riemann curvature of $F$. Then the $h$-curvature tensor $H^{\
i}_{j\ kl}$ of Berwald connection are defined by
 $$
H^{\ i}_{j\ kl}:=\frac{1}{3}\Big(\frac{\pa^2R^i_{\ l}}{\pa y^j\pa
y^k}-\frac{\pa^2R^i_{\ k}}{\pa y^j\pa y^l}\Big).
 $$
Further we define
$$H_{jk}:=H^{\ p}_{j\ kp},\ \
H_j:=\frac{1}{n-1}(nH_{0j}+H_{j0}),$$
 and then the coefficients  $K_{ij}$ of the $K$-curvature are
 define as
  \be\label{ycw111}
K_{ij}:=H_{i;j}-H_{j;i},
  \ee
where the symbol $_{;j}$ denotes the $h$-derivative of  Berwald
connection. It is shown in \cite{Ma1} that if a Finsler metric is
of scalar flag curvature $\lambda$, then we have
 \be\label{ycw112}
 H_i=(n+1)\big(\frac{1}{3}F^2\lambda_{y^i}+\lambda FF_{y^i}\big).
 \ee
In \cite{Ma1}, it proves that a two-dimensional Finsler metric $F$
is locally projectively flat if and only if $F$  is Douglasian and
the $K$-curvature vanishes $K_{12}=0$.

\begin{ex}\label{ex3}
 Let $F=c\beta+\alpha^2/\beta$, where $c$ is a constant. Define
  $$
  \alpha=\eta^{\frac{m}{m-1}}\sqrt{(y^1)^2+(y^2)^2}, \ \ \
  \beta=\eta y^1,
   $$
where $\eta=\eta(x^2)$. Then $F$ is locally projectively flat if
and only if $c\eta'''=0$.
 \end{ex}

{\it Proof :} By (\ref{ycw111}) and (\ref{ycw112}), a direct
computation gives
 $$
 K_{12}=-\frac{3}{2}c\eta'''y^1.
 $$
Now it is clear that $K_{12}=0$ if and only if $c\eta'''=0$.  \qed

\begin{ex}
 Let $F=c\beta+\alpha^4/\beta^3$, where $c\ne 0$ is a constant. In
 (\ref{yc1}) and (\ref{yc2}), put
  $$
  u=x^2,\ \ v=-x^1,\ \ B=x^1,\ \ c=1,
   $$
   and let $\alpha$ and $\beta$ be defined by (\ref{yc1}) and
   (\ref{yc2}).
 Then by Theorem \ref{th001}(ii), $F$ is a Douglas metric. However,
 $F$ is not locally projectively flat.
 \end{ex}

{\it Proof :} Similarly as in the proof to Example \ref{ex3}, we
only need to compute $K_{12}$. A direct computation gives
 $$
K_{12}=\frac{3(A_1y^1+A_2y^2)}{[4+3(x^1)^2]^5[(x^1)^2+(x^2)^2]^2},
 $$
where $A_1,A_2$ are defined by
 \beqn
 A_1:&=&d\big\{1296d^7e+e(3555+540e^2)d^5+e(720e^2+2820)d^3+e(224-960e^2)d-1280e^3\big\},\\
  A_1:&=&d\big\{-540d^8+(216e^2-2115)d^6+(720e^2-3012)d^4+(768e^2-1248)d^2+256e^2\big\},
 \eeqn
where $d:=x^1,e:=x^2$. Now it is clear that $K_{12}\ne 0$. So $F$
is not locally projectively flat.

\vspace{0.6cm}

\noindent Guojun Yang \\
Department of Mathematics \\
Sichuan University \\
Chengdu 610064, P. R. China \\
{\it  e-mail :} ygjsl2000@yahoo.com.cn

\end{document}